\documentclass{article}
\usepackage{amsthm}
\usepackage{amsmath}
\usepackage{amssymb}
\usepackage{epstopdf}
\usepackage{graphicx}

\begin{document}
\title
{Sandwiching the Riemann hypothesis}
\author{R.C. McPhedran,\\
School of Physics,\\
University of Sydney}
\maketitle
\begin{abstract}
We consider a system of three analytic functions, two of which are known to have all their zeros on the critical line $\Re (s)=\sigma=1/2$. We construct inequalities which constrain the third function, $\xi(s)$, on $\Im(s)=0$ to lie between the other two functions, in a sandwich structure. We investigate what can be said about the location of zeros and radius of convergence of expansions of  $\xi(s)$, with promising results.
\end{abstract}

\section{Introduction} 
The Riemann hypothesis  has a long history and important consequences in number theory, which are treated in numerous textbooks e.g. \cite{titheath,edw}. Riemann's hypothesis
was that all the non-trivial zeros of the Riemann zeta function $\zeta (s)$ lie on the critical line $\Re(s)=1/2$. There have been many methods developed in order to try to prove this result by analytic reasoning and by impressive computer evaluations (with the hypothesis verified numerically up to $3 \times 10^{12}$ \cite {plattrud}).

Among the many threads followed in the literature, we concentrate here on one based on sums over inverse powers of the zeros of $\zeta (s)$ on the critical line \cite{lehmer,keiper,li}.
This is chosen because it permits a combination of numerical exploration and analytic reasoning, with significant examples of both in the literature. The focus of this work will be on the
Keiper-Li theory, developed independently by J. B. Keiper and X. J. Li. This is based on expansions of the symmetrized counterpart of $\zeta (s)$, $\xi (s)$. The expansions involve two sets of constants, $a_n$ and $\lambda (n)$, the first due to Li, and the second to both authors. Li proved that the Riemann hypothesis holds if and only if  all  $\lambda (n)$ are non-negative. The author and colleagues have recently developed an expression for the positive constants $a_n$, and evaluated to good accuracy the first 4000 of them \cite{HAL, rmcp2024}.

We will combine the Keiper-Li framework with results from the literature concerning two functions, denoted here as $\xi_+(s)$ and $\xi_-(s)$. These consist of combinations of 
$\xi(s+1/2)$ and $\xi(s-1/2)$, and it is known that all their zeros lie on the critical line and are simple \cite{prt, bomblag, lagandsuz, ki}. We will consider the properties of
$\xi_+(s)$ and $\xi_-(s)$, and show that along the positive real axis $\xi (s)$ is constrained to lie between $\xi_+(s)$ and $\xi_-(s)$, in what we will call a sandwich. We will reason from the properties of the sandwich that, under a mapping of the complex variable $s=1/2+i t$ onto the unit circle, the radius of convergence of the expansion of $\xi (s)$ is unity. This result is highly significant with respect to the Riemann hypothesis.

\section{Basic equations and some results}
We will be interested here in some functions related to the Riemann $\xi$ function:
\begin{equation} 
\xi(s)=\frac{1}{2} s(s-1)\frac{\Gamma (s/2)\zeta(s)}{\pi^{s/2}}.
\label{sec2eq1}
\end{equation}
(Note that Li \cite{li} uses a definition of $\xi(s)$ which is two times larger than that used here.)
This is an entire function, with $\log \xi (s)$ having logarithmic singularities at the roots $\rho$ of $\xi(s)$ and no other singularities \cite{edw}. The particular functions of interest
are $\xi(s)$, $\xi(s+1/2)$, $\xi(s-1/2)$ and two combinations of the latter two:
\begin{equation}
\xi_{-}(s)=\xi(s+1/2)-\xi(s-1/2), ~~\xi_{+}(s)=\xi(s+1/2)+\xi(s-1/2).
\label{sec2eq2}
\end{equation}

The functions $\xi_{-}(s)$ and $\xi_{+}(s)$ are notable in that it has been proved that these functions obey the Riemann hypothesis: all their zeros $\rho$ lie on the critical line $\Re (s)=\sigma=1/2$. The functions were studied by Taylor \cite{prt}, Lagarias and Suzuki \cite{lagandsuz} and Ki\cite{ki}. The  fact that all the non-trivial zeros  of the $\xi_{-}(s)$ combination lie on the critical line was first established by P.R. Taylor, and published  
posthumously. Lagarias and Suzuki considered the $\xi_{+}(s)$ combination, and showed that all its complex zeros lie on the critical line, while Ki proved that all the complex zeros were simple. A further useful property is that the complex zeros of the $\xi_{+}(s)$ and $\xi_{-}(s)$ combinations strictly alternate on the critical line, and have the same distribution function of zeros. The common distribution function is indeed that corresponding to any prescribed argument value of $\xi(s)$ on the line $\sigma=1$.

Lagarias and Suzuki \cite{lagandsuz} in their Theorem 4 present a general result relating to functions consisting of a superposition of two identical parts, each of the parts having zeros symmetrically placed about a midline in a critical strip, and the two parts separated sufficiently so their critical strips do not overlap.  Then the symmetry of the superposition guarantees that  the moduli of the two parts can only be equal halfway between the two separated critical regions. In turn, if the superposition of the two parts is a sum with a phase factor of unit  modulus, that sum will have all its zeros on the midline of the total system. This result then includes Taylor's 
result for $\xi_{-}(s)$  as a special case, along with that for  $\xi_{+}(s)$. Related previous work is \cite{titheath,jacques3}.

We investigate the properties of the five $\xi$ functions using power series known to have coefficients all of which are known to be positive. The first has been studied intensively by
Pustyl'nikov \cite{pust,pust2} and subsequent authors. The set of coefficients $\xi_r$ occurs in the following expansion:
\begin{equation}
 \xi (s+1/2)=\sum_{r=0}^\infty \xi_r s^{2 r},
\label{sec2eq3}
\end{equation}
where $\xi_0=\xi(1/2)\approx 0.49712077818831410991$.
Pustyl'nikov proved that all the $\xi_r$ are positive and that this is a necessary condition for the Riemann hypothesis to hold. Accurate numerical techniques for the evaluation of the
$\xi_r$ have been developed in subsequent work- see for example \cite{GORZ}. A comprehensive and accurate tabulation of values due to Dr. Rick Kreminski \cite{kreminski} was previously available on the Internet, but it appears to be currently inaccessible.

The expansion (\ref{sec2eq3}) can be re-expressed to give
\begin{equation}
\xi_+(s)=2 \sum_{n=0}^\infty \xi_n s^{2 n}-\sum_{n=0}^\infty \sum_{r=n+1}^\infty \left[s\binom{2 r}{2 n+1}-\binom{2 r}{2 n}\right]\xi_r s^{2 n} ,
\label{sec2eq4}
\end{equation}
and
\begin{equation}
\xi_-(s)=\sum_{n=0}^\infty \sum_{r=n+1}^\infty \left[s\binom{2 r}{2 n+1}-\binom{2 r}{2 n}\right]\xi_r s^{2 n} .
\label{sec2eq5}
\end{equation}
These are respectively even and odd under the substitution $s\rightarrow 1-s$. Adding and subtracting (\ref{sec2eq4}) and (\ref{sec2eq5}) gives the expansion (\ref{sec2eq3})  for $\xi(s+1/2)$ and  the corresponding equation for $\xi(s-1/2)$.

The corresponding result for $2\xi(s)$ also contains both odd and even powers of $s$ and can be written
\begin{equation}
2\xi(s)=1+\sum_{l=1}^\infty {\cal E}_l s^{2 l}-\sum_{l=1}^\infty {\cal O}_l s^{2 l-1},
\label{sec2ext1a}
\end{equation}
where
\begin{equation}
{\cal E}_l =2 \sum_{r=l}^\infty{2 r\choose 2l}  \frac{\xi_r}{2^{2 r-2 l}}, ~ {\cal O}_l=2 \sum_{r=l}^\infty{2 r\choose 2l-1}  \frac{\xi_r}{2^{2 r-2 l+1}} .
\label{sec2ext1b}
\end{equation}
The first ten coefficients ${\cal E}$ of even order are
\begin{eqnarray}
0.023343864534226183135, 0.00025318173031652700506, 
1.7209870418615355778*10^{-6}, && \nonumber \\
 8.3159682500277216307*10^{-9}, 
3.0655602327633313510*10^{-11}, 9.0229664497612087603*10^{-14}, && \nonumber \\
2.1893251340686846583*10^{-16}, 
 4.4843405072454944930*10^{-19}, 
7.8974339566658717737*10^{-22}, && \nonumber \\
1.2134779622875435114*10^{-24 }. &&
\label{ecoeffxi}
\end{eqnarray}
The first ten coefficients $ {\cal O}$ of odd order are
\begin{eqnarray}
0.023095708966121033814, 0.00049798384992294867235, 
5.0502547922191741696*10^{-6}, && \nonumber \\
 3.2378414618810769603*10^{-8}, \
1.4852419214918940045*10^{-10}, 5.2238407222768796166*10^{-13}, && \nonumber  \\
1.4729420831622495667*10^{-15}, 3.4352126539793423994*10^{-18}, 
6.7821598786771781572*10^{-21}, && \nonumber \\
 1.1540437076606000624*10^{-23}.
\label{ocoeffxi}
\end{eqnarray}

The second set of power series is based around the Li coefficients $a_n$ \cite{li}. These occur in the expansion
\begin{equation}
2 \xi \left( \frac{1}{1-z}\right)=\phi(z)=1+\sum_{j=1}^\infty a_j z^j,
\label{sec2eq6}
\end{equation}
for $z$ in the unit disc. An equivalent form for $w$  to the right of the critical line is obtained using the substitution $1-z=1/w$:
\begin{equation}
2 \xi (w)=\phi\left(\frac{w-1}{w}\right)=1+\sum_{j=1}^\infty a_j \left(\frac{w-1}{w}\right)^j=\phi_1(w).
\label{sec2eq7}
\end{equation}
To the left of the critical line,  using the symmetry property $\xi(s)=\xi(1-s)$, we find:
\begin{equation}
2 \xi (w)=\phi\left(\frac{w}{w-1}\right)=1+\sum_{j=1}^\infty a_j \left(\frac{w}{w-1}\right)^j=\phi_2(w).
\label{sec2eq7a}
\end{equation}
Note that $\phi_1(1)=1$ exactly, and $\phi_2(0)=1$ exactly. 

Li's paper \cite{li} is important in that it establishes as a necessary condition for the Riemann hypothesis to hold that the logarithmic derivative of the function $\xi(1/(1-z))$ be analytic in the unit disc: see also Bombieri and Lagarias \cite{bomblag} and Coffey \cite{coffey2005}. Important related work by Keiper \cite{keiper} predated that of Li. Note that the following
expression has been established for the $a_n$ \cite{HAL, rmcp2024}:
\begin{equation}
 a_n=2\sum_{p=1}^n {\cal C}_{n,p}  \Sigma_p^\xi ,
\label{sec2eq8}
\end{equation}
where the $\Sigma_p^\xi $ are positive-valued sums involving the $\xi_r$:
\begin{equation}
 \Sigma_p^\xi=\sum_{r=1}^\infty \frac{\xi_r}{2^{2 r}} r^p,
 \label{sec2eq9}
 \end{equation}
 and $\Sigma_0^\xi =1/2-\xi_0$.
 The  $C_{n,p}$ arise in coefficients of polynomials $a_r(n)$ generated by the expansion of a quotient function:
 \begin{equation}
\left(\frac{1+w}{1-w}\right) ^r=1+\sum_{n=1}^\infty a_r(n) w^n.
\label{sec2eq10}
\end{equation}
All the  $C_{n,p}$ are non-negative, being zero if the pair $n,p$ mixes an even and an odd integer, and being positive if both are even or odd.  Note that $C_{n,n}=4^n/n!$ with $C_{n,p}$ tending to zero as $n$ increases. The  $C_{n,p}$ obey a simple recurrence relation  \cite{HAL, rmcp2024}, from which they may be determined exactly.

The equation (\ref{sec2eq8}) has been used to obtain the first 4000 $a_n$ values to 200 decimal places accuracy \cite{HAL, rmcp2024}, and they may be downloaded. The first 500
$\xi_r$ values were sufficient to supply this accuracy. Note that
\begin{equation}
a_1=\lambda_1=\sigma_1=1+\gamma/2-\log (4 \pi)/2 =8 \Sigma_1^\xi.
\label{sec2eq11}
\end{equation}
Here $a_1\approx 0.023095708966121033814$.

The ${\cal C}_{n,p}$ obey the useful relation: 
\begin{equation}
\sum_{p=1}^n {\cal C}_{n,p} =4 n.
\label{sec2eq12}
\end{equation}
Also, from equation (\ref{sec2eq9}) the $\Sigma_p^\xi $ increase monotonically with $p$. Hence, the $a_n$ for $n>1$ satisfy the bounds
\begin{equation}
n a_1<a_n<8 n \Sigma_n^\xi.
\label{sec2eq13}
\end{equation}
From \cite{HAL} $\log \Sigma_n^\xi$ has the leading terms for  $n$ large
\begin{equation}
n [\log n -\log\log n -2]
\label{sec2eq14}
\end{equation}
so the lower and upper bounds in (\ref{sec2eq13}) differ considerably in their growth rates as $n$ increases.

One important consequence of the inequalities (\ref{sec2eq13}) is that 
\begin{equation}
\sum_{n=1}^Na_n >a_1 \sum_{n=1}^N n=a_1 \frac{N(N+1)}{2},
\label{sec2eq13a}
\end{equation}
which diverges quadratically as $N\rightarrow \infty$. A second is that the series $ \phi_1(w)$ must diverge if $|(w-1)/w|>1$, i.e. if $Re(w)<1/2$ and $|z|>1$, while the series $ \phi_2(w)$ must diverge if $|w/(w-1)|>1$,   i.e. if $Re(w)>1/2$ and $|z|<1$.
\section{The first two sandwiches}
The equation (\ref{sec2eq7}) can be used to construct inequalities and monotonicity properties among the key  functions introduced in the previous section. Consider the difference
$\xi(u+a)-\xi(u+b)$ where $u$,$a$ and $b$ are real with $a>b$ and $u+a$ and $u+b$ both exceeding 1/2 :
\begin{equation}
2 \xi (u+a)-2 \xi (u+b)=\sum_{j=1}^\infty a_j \left[\left(1-\frac{1}{u+a}\right)^j-\left(1-\frac{1}{u+b}\right)^j\right]>0,
\label{sec3eq1}
\end{equation}
since every element of the summand is positive. Hence, $\xi(u+1/2)>\xi(u)$ and  $\xi(u)>\xi(u-1/2)$. For $\xi(u-1/2)$ we have 
\begin{equation}
2 \xi (u-1/2)=1+\sum_{j=1}^\infty a_j \left(1-\frac{1}{u-1/2}\right)^j> 1,~{\rm for}~u>3/2.
\label{sec3eq2}
\end{equation}
 The full set of inequalities for the first three functions of interest is then
\begin{equation}
\xi\left(u+\frac{1}{2}\right)>\xi(u)>\xi\left(u-\frac{1}{2}\right) >0.
\label{sec3eq3}
\end{equation}
Each of the three functions increases monotonically as the argument $u$ increases, and also if the full sum is replaced by a partial sum, the partial sums increase monotonically  with the upper limit on $j$.

We now consider mappings  which move the location of lines along which zeros are located onto circles in the complex plane.
For the case of  $\xi(s)$ a convenient  mapping from the  critical line $\Re (s)=\sigma=1/2$ onto the unit circle has already been given:
\begin{equation}
w=u+i v=1-\frac{1}{s}=\frac{s-1}{s}.
\label{tr1}
\end{equation}
The inverse transformation is
\begin{equation}
s=\frac{1}{1-w}.
\label{tr2}
\end{equation}

For the function $\xi\left( s+1/2\right)$ the corresponding forward transformation is
\begin{equation}
w_h=1-\frac{1}{s+1/2}=\frac{s-1/2}{s+1/2}.
\label{tr3}
\end{equation}
Its  inverse transformation is
\begin{equation}
s=-\frac{1}{2}+\frac{1}{1-w_h}.
\label{tr4}
\end{equation}
The known zeros of $\xi\left( s+1/2\right)$ lie on $\sigma=0$ and are mapped onto the unit circle in the plane of complex $w_h$.

For the function $\xi\left( s-1/2\right)$ the corresponding forward transformation is
\begin{equation}
w_m=1-\frac{1}{s-1/2}=\frac{s-3/2}{s-1/2}.
\label{tr5}
\end{equation}
Its  inverse transformation is
\begin{equation}
s=\frac{1}{2}+\frac{1}{1-w_m}.
\label{tr6}
\end{equation}
The known zeros of $\xi\left( s-1/2\right)$ lie on $\sigma=1$ and are mapped onto the unit circle in the plane of complex $w_m$.

We next connect $w_h$ and $w_m$  to the complex variable $w$. Eliminating $s$ between (\ref{tr1}) and (\ref{tr4}) we find
\begin{equation}
w_h=\frac{1+w}{3-w}, ~~w=\frac{3 w_h-1}{w_h+1}.
\label{tr7}
\end{equation}
The equation for the fixed point of this transformation is
\begin{equation}
w (3-w)=1+w, ~{\rm or}~ (w-1)^2=0.
\label{tr8}
\end{equation}
The fixed point is then of second order at $w=w_h=1$. The corresponding equations relating to $w_m$ are
\begin{equation}
w_m=\frac{3 w-1}{w+1}, ~~w=\frac{1+w_m}{3-w_m},
\label{tr9}
\end{equation}
while again the fixed point for the transformation yields $w=1=w_m$ being of second order.

\begin{figure}[tbh]
\includegraphics[width=7.5 cm]{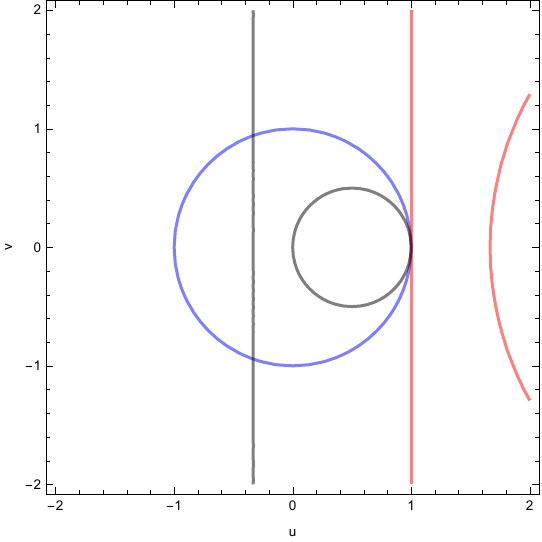}
\caption{ The transformations pertaining to $w$ (blue circle), $w_m$ (black circle and black line) and $w_h$ (red lines).}
\label{fig-regions}
\end{figure}

Figure \ref{fig-regions} illustrates curves of constant modulus in the complex $w$ plane pertaining to these discussions. The black unit circle corresponds to $|w|=1$, and thus to the mapping of the critical line through equation (\ref{tr1}). The black circle centred on $w=1/2$ of radius $1/2$ corresponds to the constraint $|w_m|=1$, and the mapping of the line $\sigma=1$, with the black vertical line being $|w_m|=3$, $u=-1/3$.  The red lines are for $|w_h|=1$ and $|w_h|=2$. The two fixed points occur where the two circles touch, with the line $|w_h|=1$ being tangent to both.

\begin{figure}[tbh]
\includegraphics[width=12cm]{"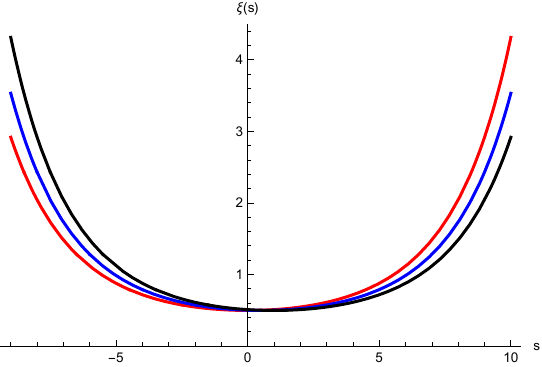"}
\caption{The behaviour of the three functions $\xi(s+1/2)$ (red), $\xi(s)$ (blue) and $\xi(s-1/2)$ (black) along part of the real axis is shown.}
\label{figxisimp}
\end{figure}

The behaviour of the three functions $\xi(s+1/2)$, $\xi(s)$ and $\xi(s-1/2)$ along a part of the real axis is shown in Fig. \ref{figxisimp}, and is in accord with the inequalities (\ref{sec3eq3}) . The figure shows the first of the two sandwiches we will study. Our aim will be to understand the location of the singularities of the three functions, as reflected in their power series and their behaviour
along the real axis.

Titchmarsh \cite{titfn} gives two useful results concerning the radius of the circle of convergence of power series with coefficients $a_n$ in the variable $z$ and their behaviour along the real axis. For convenience the radius of the circle of convergence is normalised to be unity. The first is that if $a_n\geq 0$ for all $n$, then $z=1$ is a singular point. The second is that if $a_n$ is real for all values of $n$ and $\sum a_n$ is properly divergent, then $z=1$ is a singular point. We note that both criteria apply to power series with the coefficients which are the Li $a_n$ (see equation (\ref{sec2eq13a})). If the coefficients are the $\xi_n$, the first will apply.

We thus deduce that the circle of convergence of $\phi(z)$ cuts the positive real axis at a singular point of its power series (\ref{sec2eq6}). The same conclusion applies to the circles of convergence pertaining to $\phi(z_h)$ where $z_h=z+1/2$ and to $\phi(z_m)$ where $z_m=z-1/2$.

\begin{figure}[tbh]
\includegraphics[width=12cm]{"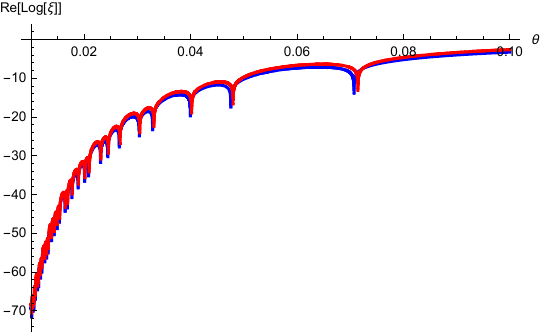"}
 \includegraphics[width=12cm]{"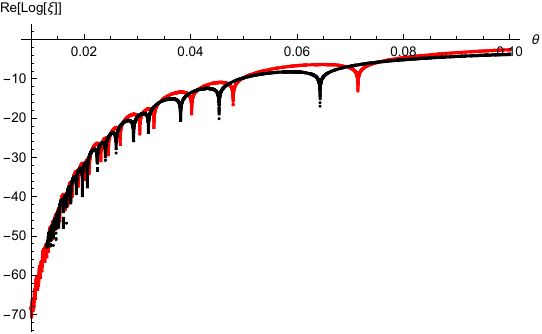"}
\caption{The behaviour of the three functions $\xi_+(s)$ (red), $\xi(s)$ (blue) and $\xi_-(s)$ (black) along part of unit circle is shown.}
\label{figxiscirc}
\end{figure}

The behaviour of the three functions $\xi_+(s)$, $\xi(s)$ and $\xi_-(s)$ along an arc of the unit circle is shown in Fig. \ref{figxiscirc}. Specifically, the function plotted is $\log\xi[1/(1-\exp(i \theta)]$, where $\theta$ denotes the angular position on the unit circle. The dips in the graph indicate zeros of the function. For $\xi_+$ they occur where 
$\arg[\xi(1 + (i/2) \cot (\theta/2))]$ is an odd multiple of $\pi/2$, and for $\xi_+$ where it is an even multiple of $\pi/2$. Note that zeros of  $\xi_+(s)$ occur closer to those of $\xi(s)$
than is the case for zeros of  $\xi_-(s)$, and that  zeros of  $\xi_+(s)$ strictly alternate with those of $\xi_-(s)$. (This is a consequence of the monotonic variation of $\arg[\xi(1 + (i/2) t]$ with $t$.)

\begin{figure}[tbh]
\includegraphics[width=12cm]{"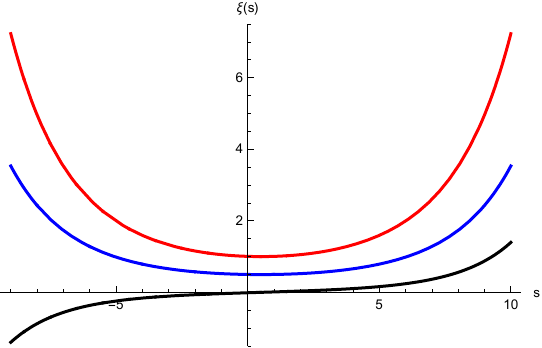"}
\caption{The behaviour of the three functions $\xi_+(s)$ (red), $\xi(s)$ (blue) and $\xi_-(s)$ (black) along part of the real axis is shown.}
\label{figxis}
\end{figure}

The behaviour of  $\xi_+(s)$, $\xi(s)$ and $\xi_-(s)$ along a part of the real axis is shown in Fig. \ref{figxis}. The figure gives evidence for the following inequalities:
\begin{equation}
\xi_+(\sigma)>\xi(\sigma)>\xi_-(\sigma)>0~{\rm for}~ \sigma\geq 1.
\label{sec3eq4}
\end{equation}
These inequalities follow immediately from (\ref{sec3eq3}), with the exception of $\xi(\sigma)>\xi_-(\sigma)$. This last follows from the asymptotic expression:
\begin{eqnarray}
2\xi(\sigma)-2[\xi(\sigma+1/2)-\xi(\sigma-1/2)]&=&1+\sum_j a_j \left(1-\frac{1}{\sigma}\right)^j \nonumber \\
&& \left\{ 1+\exp\left[ j\left( \log \left(1-\frac{1}{(\sigma-1/2)}\right)- \log \left(1-\frac{1}{\sigma}\right)\right) \right]  \right. \nonumber \\
&& \left.- \exp\left[j \left(\log \left(1-\frac{1}{(\sigma+1/2)}\right)- \log \left(1-\frac{1}{\sigma}\right)\right) \right] \right\} \nonumber \\
&&
\label{sec3eq5}
\end{eqnarray}
giving for $j<<\sigma^2$:
\begin{eqnarray}
 2\xi(\sigma)-2[\xi(\sigma+1/2)-\xi(\sigma-1/2)]=1+\sum_j a_j \left(1-\frac{1}{\sigma}\right)^j  
\left[ 1-\frac{j}{\sigma ^2}-\frac{j}{\sigma ^3}+O\left(\left(\frac{1}{\sigma }\right)^4\right)\right] &&
\label{sec3eq6}
\end{eqnarray}
This is  increasing with $\sigma$, diverging as $\sigma\rightarrow \infty$ by virtue of the increasing values of $a_j$. The functions $\xi_+(\sigma)$, $\xi(\sigma)$ and $\xi_-(\sigma)$ also increase monotonically with $\sigma>1$ and are all positive for $\sigma>1/2$ . They must all diverge for $\sigma\rightarrow \infty$. 

\section{The third sandwich}
\begin{figure}[tbh]
\includegraphics[width=12cm]{"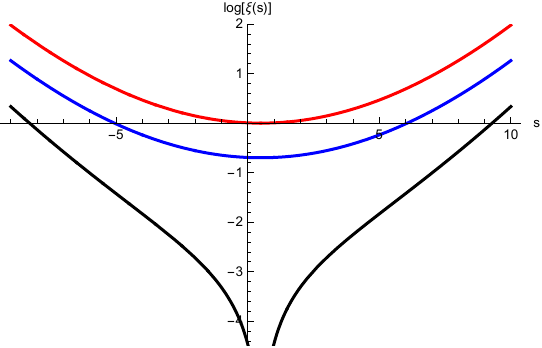"}
\caption{The behaviour of the three functions $\log(\xi_+(s))$ (red), $\log(\xi(s))$ (blue) and $\Re \log(\xi_-(s))$ (black) along part of the real axis is shown.}
\label{figlogxis}
\end{figure}

We next consider the logarithms of the three functions  $\xi_+(s)$, $\xi(s)$ and $\xi_-(s)$, as shown in Fig. \ref{figlogxis}. The curves there illustrate a similar behaviour to previous figures for  $\xi_+(s)$ and $\xi(s)$. For $\xi_-(s)$ there is a notable influence of the first order zero at $s=1/2$. The monotonic behaviour of the logarithm enables us to go from
the inequalities (\ref{sec3eq4}) to 
\begin{equation}
\log \xi_+(\sigma)>\log \xi(\sigma)>\log \xi_-(\sigma)>0~{\rm for}~ \sigma\geq 1.
\label{sec4eq1}
\end{equation} 
The monotonic increase of the three functions of (\ref{sec4eq1}) with $\sigma$ is also preserved.

A power series for $\log \xi(s)$ with coefficients  depending on the equation (\ref{sec2ext1a}) can be constructed in two ways. The first is to use Mathematica or a similar symbolic 
package. The second is to differentiate the series for $\xi(s)$ and divide the result by  $\xi(s)$, before integrating the result for the quotient. The result is the same, to within  slight numerical differences:
\begin{eqnarray}
\log(\xi (s))=  -0.023095708966121033814 s + 0.023077158647902301379 s^2 + 
 0.0000370527438173686409 s^3 && \nonumber \\
 - 0.00001840680531542237958 s^4 - 
 1.43018671152521547*10^{-7 }s^5 + 4.6906069489794377*10^{-8} s^6  && \nonumber \\
+  6.534558785292532*10^{-10} s^7 - 1.5860851386884509*10^{-10} s^8 - 
 3.141596838950986*10^{-12} s^9 &&  \nonumber \\
  + 5.99771484715187*10^{-13} s^{10} + 
 1.53698978205383*10^{-14} s^{11} - 2.39484594174383*10^{-15} s^{12}  && \nonumber \\
 - 7.5638205856654*10^{-17} s^{13 }+ 9.851878115521*10^{-18} s^{14} + 
 3.727558439514*10^{-19} s^{15}  && \nonumber \\
 - 4.123828024852*10^{-20} s^16 - 
 1.836227510375*10^{-21} s^{17} + 1.743694268780*10^{-22} s^{18}  && \nonumber \\
  +  9.03450466141*10^{-24} s^{19} - 7.4119966493*10^{-25} s^{20}
 \label{serlxi}
 \end{eqnarray}
 This is a power series with coefficients of mixed sign, unlike that for $\xi(s)$. It has been obtained using Mathematica, by the first method, with 20 decimal places requested as the accuracy goal. Note that higher order coefficients have fewer decimals 
than specified in the accuracy goal.

\begin{figure}[tbh]
\includegraphics[width=8cm]{"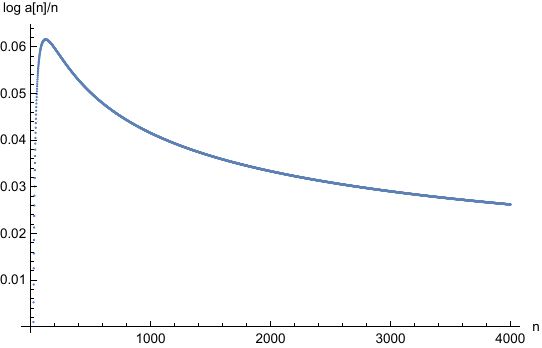"}~\includegraphics[width=8cm]{"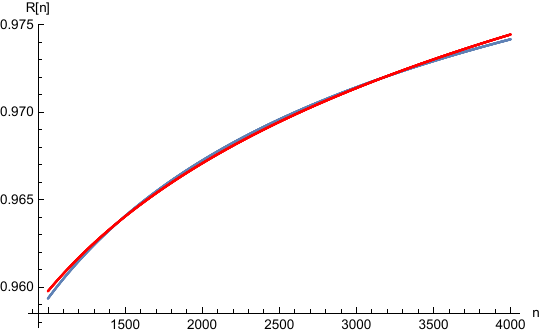"}
\caption{(Left) The values of $\log a_n/n$ for $n$ ranging up to 4000. (Right) The  $R_n$ values (blue) and the fit function (\ref{sec4eq3}) (red) for $n$ ranging from 1000 to 4000. }
\label{figlogxisa}
\end{figure}

The main focus of this section will be a discussion of the radius of convergence of the expansion (\ref{sec2eq6}), and its equivalent 
for $\log[2 \xi ( 1/(1-z))]$. From \cite{titfn}, the radius of convergence associated with the series in (\ref{sec2eq6}) is 
\begin{equation}
R=\lim_{n\rightarrow \infty} R_n=\lim_{n\rightarrow \infty} \exp\left[\frac{-\log  a_n}{n}\right].
\label{sec4eq2}
\end{equation}
 In Fig.  \ref{figlogxisa} we show  the variation of  $\log a_n/n$ with $n$,  with a rapid rise for  values of $n$ up to around 300 being succeeded by a slow fall off for higher values.
Also shown are the associated values of $R_n$ from equation (\ref{sec4eq2})
and a numerical fit for the range of $n$ from 1000 to 4000:
\begin{equation}
\frac{\log a_n}{n}\approx 0.11652745618 -0.01092578334 \log n .
\label{sec4eq3}
\end{equation}
The simple fit function gives a reasonably accurate representation of the variation of $\log a_n/n$ in the range shown, and confirms that its leading varying term is logarithmic. Values of $a_n$ for $n$ much larger than 4000 would be needed to get additional terms in the expression for $\log a_n/n$  and thus for $R_n$.

\begin{figure}[tbh]
\includegraphics[width=8cm]{"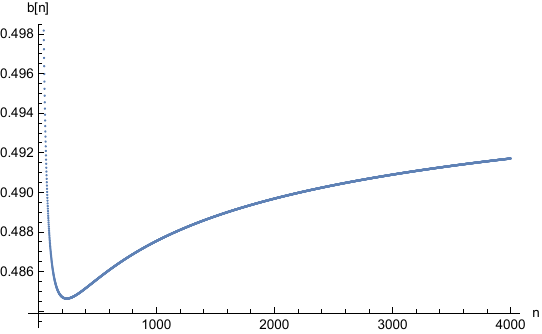"}
\caption{The values of the singularity estimate $b_n$ as a function of $n$. }
\label{Titbnplt}
\end{figure}

Titchmarsh \cite{titfn} gives a necessary and sufficient condition that a power series with coefficients $a_n$ (known here to be non-negative) should have a singularity at the point $z=1$ lying on its circle of convergence. Let $b_n$ denote the quantities 
\begin{equation}
b_n=\sum_{m=0}^n {n \choose m} a_m.
\label{sec4eq4}
\end{equation}
Then the criterion for a singularity to occur at $z=1$ is
\begin{equation}
\lim_{n\rightarrow \infty} b_n^{-\frac{1}{n}}=\frac{1}{2}.
\label{sec4eq5}
\end{equation}
The application of this criterion to the power series with coefficients $a_n$ is shown in Fig. \ref{Titbnplt}. While both the radius of convergence estimate and the singularity estimate
converge slowly with increasing $n$, they are confirmatory of there being a singularity at $z=1$ lying on the unit circle of convergence.

Of more importance to the following discussion than the series for   $\log \xi(s)$  are those for  $2 \xi \left( \frac{1}{1-z}\right)=\phi(z)$ and $\phi'(z)/\phi(z)$. For the last, Li states that a necessary and sufficient condition for the nontrivial zeros of $\zeta(s)$ to lie on the critical line is that $\phi'(z)/\phi(z)$ is analytic in the unit disc. Furthermore, given that the series for $\phi'(z)/\phi(z)$ around $z=0$ is
\begin{equation}
\frac{\phi'(z)}{\phi(z)}=\sum_{n=0}^\infty \lambda_{n+1} z^n,
\label{Lilamb}
\end{equation}
the necessary and sufficient condition is equivalent to the requirement that $\lambda_n\geq 0$ for every positive $n$. Keiper had previously shown that non-negativity of  the   
$\lambda_n$ was a necessary condition for the Riemann hypothesis, and had calculated the $\lambda_n$ to high accuracy up to $n=4000$. Note that Keiper based his definition of the $\lambda_n$ around the expansion of $\log \phi (z)$, resulting in a factor of $n$ between his $\lambda_n^K$ and those of Li ($\lambda_n^L$).

One way of going from equation (\ref{sec2eq6}) to the required series for $\phi'(z)/\phi(z)$ is to form the series for $1/\phi(z)$:
\begin{equation}
\frac{1}{\phi(z)}=1+\sum_{j=1}^\infty A_j z^j,
\label{myAexp1}
\end{equation}
via the recurrence relation
\begin{equation}
A_1=-a_1,~ A_j=-a_j-\sum_{p=1}^{j-1} a_p A_{j-p}.
\label{myAexp2}
\end{equation}
The required series then comes from the product:
\begin{equation}
[\sum_{j=1}^J j a_j z^{j-1}] \times [1+\sum_{j=1}^J A_j z^j],
\label{myAexp3}
\end{equation}
where $J$ specifies the number of values for $\lambda$'s occurring 
in the series.

\begin{figure}[tbh]
\includegraphics[width=8cm]{"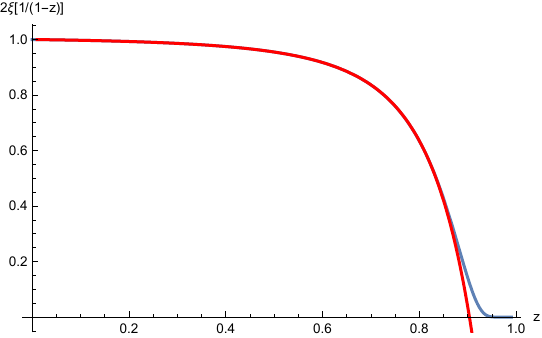"}
\caption{The function $1/\phi(z)$ as a function of $z$ is compared with the series (\ref{resltrecip20}). }
\label{ratphiponphiplt}
\end{figure}

The series for $1/\phi(z)$ is to order 20:
\begin{eqnarray}
  1-0.0230957089661210338143102479065 z-0.0459061617276994534365358813998  z^2 & & \nonumber \\
   -0.0681486316594069122599620826562 z^3-0.089545433048398089752378459144 z^4 & & \nonumber \\
   -0.109826396050622458253486711461 z^5   -0.128731287368373527211506921969 z^6 & & \nonumber \\
   -0.146012159089427308368403139905 z^7-0.161435608581258419635214039431 z^8 & & \nonumber \\
   -0.174784932895325800969831596126 z^9-0.185862161803051549940505866001 z^{10}  & & \nonumber \\
   -0.194489954353028235590838150534 z^{11}  -0.200513344703746739263814823922 z^{12}  & & \nonumber \\
   -0.203801323942739062262628361287 z^{13}-0.20424824564594307240714200573 z^{14} & & \nonumber \\
   -0.20177504405427703081208059797 z^{15}-0.19633025494134585958756091470z^{16}  & & \nonumber \\
   -0.18789083050993567984657371103 z^{17} -0.17646274097813564381563835719  z^{18}   & & \nonumber \\
    -0.16208135689084447159310105477 z^{19}-0.14481160761104406228773575035 z^{20} & & \nonumber \\
 +O\left(z^{21}\right) . & & 
   \label{resltrecip20}
   \end{eqnarray}
  Here, all coefficients of $z$ are negative. This results in a monotonic decreasing function, as shown in Fig. \ref{ratphiponphiplt}. Note that the truncated series (\ref{resltrecip20}) works well until $z$ approaches the exponential region near unity.
   
  The series for $\phi'(z)/\phi(z)$ is to order 19: 
  \begin{eqnarray}
\frac{\phi'(z)}{\phi(z)}=\sum_{n=0}^\infty \lambda_{n+1}^L   z^n \hspace{5 cm}& & \nonumber \\
 = 0.0230957089661210338143102479065+0.0923457352280466703857284861921 z & & \nonumber \\
  +0.207638920554324803791492046618 z^2+0.368790479492241638590511489638 z^3 & & \nonumber \\
 +0.575542714461177452431106405493 z^4+0.82756601228237929742500282202  z^5 & & \nonumber \\
  +1.12446011757095949058282010802 z^6+1.46575567714706063265551454198  z^7 & & \nonumber \\  
   +1.85091604838253415532604486792 z^8+2.27933936319315774369303405737  z^9 & & \nonumber \\
  +2.75036083822019606035454709285 z^{10}+3.26325532062461984807908598991  z^{11} & & \nonumber \\
  +3.81724005784794598710436795129 z^{12}+4.4114776786805985120806412969  z^{13} & & \nonumber \\
  +5.0450793720267934585351114375 z^{14}+5.7171082488687926394190666698 z^{15} & & \nonumber \\
   +6.4265828721172029011455409609 z^{16}+7.1724809382917229592529707263  z^{17} & & \nonumber \\
  +7.9537430943119003048082250779 z^{18}+8.7692768720932151295994613534   z^{19}& & \nonumber \\
 +O\left(z^{20}\right)
   \label{resltphirat20}
   \end{eqnarray}
   All coefficients of $z$ being positive, this is a monotonically increasing function, positive within the radius of convergence of the series.  We noted at the beginning of this section
   that the monotonic increase of the functions $\xi_+(s)$, $\xi(s)$ and $\xi_-(s)$ for $s$ real and larger than unity was preserved for  $\log \xi_+(s)$, $\log \xi(s)$ and $\log \xi_-(s)$,
   whose series incorporate the same coefficients divided by their power as do their derivatives.

More definite information comes from the knowledge that $\xi_+(s)$ and $\xi_-(s)$ have all their zeros on the critical line, mapped by $s\rightarrow 1/(1-z)$ onto the boundary of the unit circle. The point $z=1$ is then a limit point for a sequence of zeros (see Fig. \ref{figxiscirc}), and is thus an essential singularity of the logarithms of these functions.  For each
of $\log[\xi_+(\sigma)]$ and $\log[\xi_-(\sigma)]$, the function tends to infinity as $\sigma\rightarrow 1$. By virtue of the inequalities (\ref{sec4eq1}), we know then that 
$\log[\xi(\sigma)]$ tends to infinity as $\sigma\rightarrow 1$ (in keeping with  equation (\ref{sec2eq6})).

Consider what would be the situation if the Riemann hypothesis failed. There would then be at least one zero lying properly inside the unit circle, and a singularity of
 $\log[\phi(z_*)]$
at a point $z_*$ with $|z_*|<1$ (with naturally a corresponding singularity occurring outside the unit circle). Take $\rho_*$ to be the minimum of all such $|z_*|$. Then the circle of radius $\rho_*$  centred on $z=0$ is the circle of convergence of the logarithm of the  expansion in (\ref{sec2eq6}). This means that the series diverges on any ray as the modulus of $z$ approaches $\rho_*$. Applying this to the real axis, the  modulus of $\log[\phi(\sigma)]$ has to diverge to positive infinity at $\sigma=\rho_*<1$, and to exceed the modulus of $\log[\xi_+(\sigma)]$ there,  in contradiction with  (\ref{sec4eq1}).

Another way of looking at this situation is that the sandwiching shows that  $\log[\phi(\sigma)]$ cannot diverge before $\log[\xi_+(\sigma)]$ diverges. It also cannot diverge after $\log[\xi_-(\sigma)] $ diverges.

Both arguments  lead to the conclusion that indeed the radius of convergence of the function $\log[\phi(\sigma)]$ is unity, and that $\sigma=1$ is an essential singularity.
\section{Further investigations}
\begin{figure}[tbh]
\includegraphics[width=8cm]{"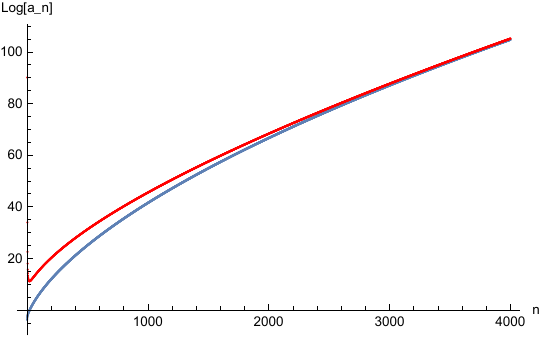"}
\caption{Numerical values of $\log (a_n)$ (blue) are compared with the asymptotic formula (\ref{finveq1}) (red). }
\label{asymloganfig}
\end{figure}

An important element of the results here is the establishment of properties of the Li coefficients $a_n$. We have given bounds on them in equation (\ref{sec2eq13}), and one numerical approximation formula in (\ref{sec4eq3}). We have found the following simple asymptotic formula for the $\log (a_n)$, which may form the basis for a more 
complete asymptotic  treatment:
\begin{equation}
\log (a_n) \sim \frac{15 n}{\log (n)^3}.
\label{finveq1}
\end{equation}
A comparison of the numerical results up to $n=4000$ with this formula is shown in Fig. \ref{asymloganfig}.

\begin{figure}[tbh]
\includegraphics[width=8cm]{"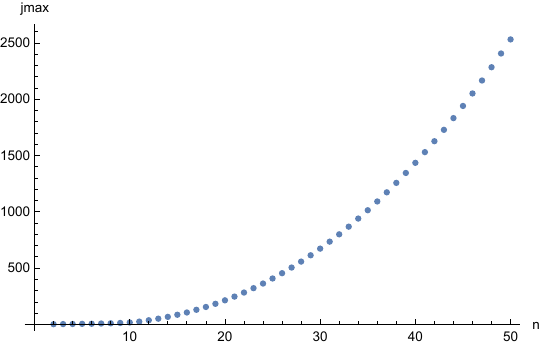"}
\caption{The value $j_m$ of $j$ which maximizes the summand in equation (\ref{sec2eq7}) is given for $n\le 50$.  }
\label{pltjm}
\end{figure}

The need for an extension of the compilation of $a_n$ values beyond $n=4000$  to give a more solid knowledge of their asymptotic behaviour is complemented by the study of the $j$ summand in the representation    (\ref{sec2eq7}) of $2\xi(s)$. The data shown in Fig. \ref{pltjm} for the location of the maximum of the summand indicates a rapid increase in the number of coefficients necessary with $n$ (the increase being slightly more than as $n^2$).

Note that the extension will require a good knowledge of the asymptotics of the coefficients ${\cal C}_{n,p}$, which is aided by a continuum treatment of the exact recurrence relation
 \cite{HAL, rmcp2024}. The continuum approach requires both $n$ and $p$ to be large, and replaces differences of ${\cal C}_{n,p}$ values by partial derivatives:
 \begin{equation}
 {\cal C}_{n-1,p-1}={\cal C}_{n,p}-\frac{\partial  {\cal C}_{n,p}}{\partial n}-\frac{\partial  {\cal C}_{n,p}}{\partial p},~ {\cal C}_{n-2,p}={\cal C}_{n,p}-2\frac{\partial  {\cal C}_{n,p}}{\partial n},
 \label{cont1}
 \end{equation}
 to first order. The recurrence relation
 \begin{equation}
 {\cal C}_{n,p}=\frac{4}{n}  {\cal C}_{n-1,p-1}+\frac{(n-2)}{n} {\cal C}_{n-2,p}
  \label{cont2}
 \end{equation}
 then reduces to the following first order identity:
 \begin{equation}
 \frac{2}{n}  {\cal C}_{n,p}=\frac{4}{n} \frac{\partial  {\cal C}_{n,p}}{\partial p}+2 \frac{\partial  {\cal C}_{n,p}}{\partial n}.
 \label{cont3}
 \end{equation}
 This has the exact solution:
 \begin{equation}
  {\cal C}_{n,p}=n {\cal F}(p-2 \log n),
 \label{cont4}
 \end{equation}
 where ${\cal F}$ is a positive-valued function with appropriate properties, for example:
 \begin{equation}
 \sum_{p=1}^n {\cal F}(p-2 \log n)=4, ~{\cal F}(n-2 \log n)=\frac{4^n}{n! n}.
 \label{cont5}
 \end{equation}
 The representation (\ref{cont4}) gives a valuable insight into how ${\cal C}_{n,p}$ depends on its two integer variables.

\begin{figure}[tbh]
\includegraphics[width=8cm]{"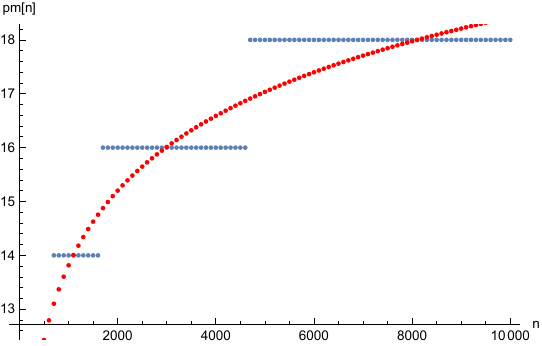"}~\includegraphics[width=8cm]{"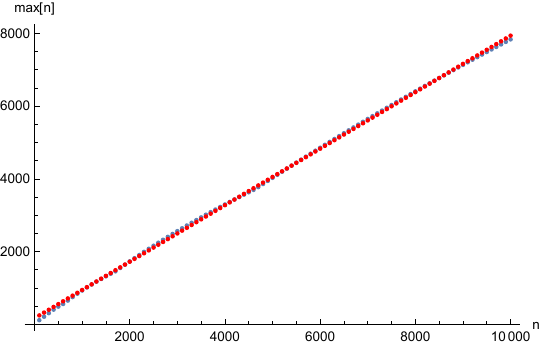"}
\caption{(left)The value $p_m$ of $p$ (blue dots) which maximizes the summand in equation (\ref{sec2eq8}) is given for $n$ in the range one to ten thousand, with the red dots
corresponding to $2\log n$. (right) The corresponding maximum values, with the red dots corresponding to the empirical formula .  }
\label{Cnppeaks}
\end{figure}

Figure  \ref{Cnppeaks} shows the variation with $n$ of the peak position $p_m$ and the corresponding peak value of the coefficients   ${\cal C}_{n,p}$, for $n$ ranging between 1000 and 10000. At left, the position of the peak is compared with the natural estimate from (\ref{cont4}). The comparison is by no means conclusive, given the logarithmic form of the estimate,
and the fact that only points from even $n$ are given, so that $\log (n)$ is required to jump by two to move from one cluster of points to the next. The value of ${\cal C}_{n,p_m}$ 
is given approximately by the empirical fit:
\begin{equation}
{\cal C}_{n,p_m}=0.78237057 n+151.978136
\label{pCnlfit}
\end{equation}

\begin{figure}[tbh]
\includegraphics[width=8cm]{"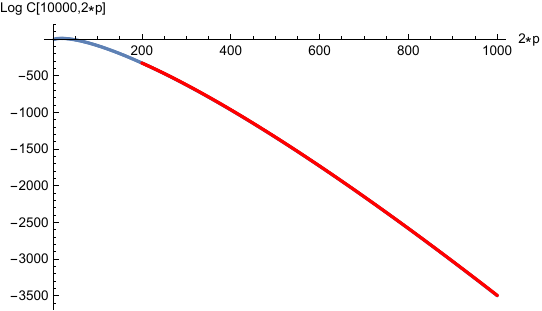"}
\caption{The fit to ${\cal C}_{n,p}$ for $n=10000$ with parameters $a,b,c$ as in the text (red points) is compared with the exact values (blue points). }
\label{Cnp10kfit}
\end{figure}

 We have for $n>>1$ the representation of ${\cal C}(n,n)/n$ or ${\cal F}(n-2 \log n)$:
 \begin{equation}
 \frac{4^n}{n! n}=\exp \left[-n \log n+(1+\log 4)n-\frac{3}{2} \log n-\frac{1}{2}\log(2 \pi)-\frac{1}{12 n}+\ldots\right].
 \label{cont6}
 \end{equation}
 We have investigated fits to  the more general case  of $\log {\cal C}_{n,p}$ for $n$ large using the following form based on (\ref{cont4}:
 \begin{equation}
 \log {\cal C}_{n,p} \sim a (p-\log n) \log (p-\log n)+b (p-\log n)+c +\log n,
 \label{cont7}
 \end{equation}
 where $a,b,c$ are the fit parameters estimated numerically.  The results for $n=10000$ are: $a = -1.188831$, $b =- 4.685604$ and 
$c = -64.9176957$, with the fit being based on values of $p$ between 200 and 1000.The fit as shown in Fig. \ref{Cnp10kfit} is good, with a maximum difference of 4.7 between the two. For smaller values of $n$, the fit parameter $a$ becomes more negative, with a slow variation approximated by $a\sim-1-2/log(n)$. The values of this parameter depend somewhat on the range of $p$ chosen- of course, the value we expect from equation (\ref{cont6}) for $p=n$  is -1. 

We next consider the numerics of the expression (\ref{sec2eq8}) for $a_n$. This contains a sum over  the product ${\cal C}_{n,p} \Sigma_p^\xi$, and data on the behaviour of this product is given in Table \ref{tabmaxansum}. Specifically, for selected values of $n$, the value of $p$, $p_a$, for which the logarithm of the summand is maximal is given, along with the value of the log summand and the value of $\log \Sigma_{p_a}^\xi$. As well, the value of the numerical derivative of $\log \Sigma_{p}^\xi$ with respect to $p$ at $p_a$ is given.
Note that $p_a$ increases with $n$ much more rapidly than the maximum location of $\log {\cal C}_{n,p} \sim 2\log n$ does. The criterion for $p_a$ is that the derivative
of the sum of these two logarithms changes sign around $p_a$. The numerical derivatives are much smaller than the values of $\log \Sigma_{p}^\xi$, and increase monotonically with
$n$ in the range shown.

\begin{table}
 \begin{tabular}{|c|c|c|c|c|}
 \hline
 $n$ & $p_a$ & $\log{\cal C}_{n, p_a} \Sigma_{p_a}$ & $\log \Sigma_{p_a}$ & N.D. $\log \Sigma_{p_a}$ \\
  \hline 
 1000 & 126 & 37.7565393217774291963&231.6146084& 2.7093408 \\
 2000 & 202 & 62.6228448499282972492&453.9908012&3.1164868 \\
 3000 & 266 & 82.8458869226227176219&661.3015708& 3.354897 \\
 4000 & 324 & 100.4888971449638584075&860.9287170& 3.526183 \\
 5000 & 376 & 116.4105045083518446758& 1047.674694& 3.655667\\
 6000 & 424 & 131.0719703207814364204&1225.655913& 3.760324 \\
 7000 & 470 & 144.756437991532648832&1400.686892& 3.850145 \\
 8000 & 514 & 157.6523866761645515334&1571.798115& 3.928260 \\
 9000 & 556 & 169.8935478487459416390& 1738.210338& 3.996876\\
 10000 & 596 & 181.5792684592106359788& 1899.283628& 4.057610\\
 11000 & 634 & 192.7855316004084369660&2054.483559& 4.111679 \\
 12000 & 672 & 203.5715544749158515816&2211.679196& 4.162628 \\
 13000 & 708 & 213.9860487621971469130& 2362.340635&4.208328 \\
 14000 & 742 & 224.0681471182688483013&2506.107283& 4.249424 \\
 15000 & 776 & 233.8503214579580060438&2651.240608& 4.288695 \\
 16000 & 810 & 243.3604498900790702158& 2797.681205& 4.326298\\
 17000 & 842 & 252.6218306046395566129& 2936.653027& 4.360286\\
 18000 & 874 & 261.6549108802165656541&3076.692759&4.393019 \\
 19000 & 906 & 270.4772209710938982525&3217.761623&4.424585\\
 20000 & 936 & 279.1040821154458906632& 3350.916248& 4.453190 \\
 \hline
\end{tabular}
\caption{The peak value of $p$, $p_a$, and of the logarithm of the  summand in  (\ref{sec2eq8}) for various values of $n$, together with the 
 value of $\log \Sigma_{p_a}$ and of its numerical derivative with respect to $p$ .}
\label{tabmaxansum}
\end{table}

We can use the recurrence relation (\ref{cont2}) for ${\cal C}_{n,p}$ to establish an exact recurrence relation for the $a_n$.  Starting from the representation( \ref{sec2eq8}) and using
(\ref{cont2}) we obtain:
\begin{equation}
a_n=\frac{4}{n} a_{n-1}+\frac{(n-2)}{n} a_{n-2}+\frac{8}{n} \sum_{p=2}^n {\cal C}_{n-1,p-1}(\Sigma^\xi_p-\Sigma^\xi_{p-1}).
\label{cont8}
\end{equation}
Since the $\Sigma^\xi_p$ increase monotonically with $p$, (\ref{cont8}) gives rise to the inequality
\begin{equation}
a_n>\frac{4}{n} a_{n-1}+\frac{(n-2)}{n} a_{n-2},
\label{cont9}
\end{equation}
and a weaker alternative:
\begin{equation}
a_n>\frac{(n+2)}{n} \min (a_{n-1}, a_{n-2}).
\label{cont10}
\end{equation}
Now, the tabulation shows the first 4000 values of $a_n$ increase monotonically with $n$. Assuming monotonicity for larger $n$ values, we can investigate an alternative to  the lower bound on $a_n$ in 
equation (\ref{sec2eq13}):
 \begin{equation}
a_n>\frac{(n+2)}{n}   a_{n-2}).
\label{cont11}
\end{equation}
For $n$ even,  (\ref{cont11}) leads to a product form resulting in :
\begin{equation}
a_{2 n}>\frac{( n + 1) }{2}  a_{2},
\label{cont12}
\end{equation}
and for $n$ odd it gives :
\begin{equation}
a_{2 n-1}> \frac{(2 n+1)}{3} a_{1}.
\label{cont13}
\end{equation}
The first of these is better as a lower bound than $n a_1$, while the second is slightly worse.

Further progress in these lines of investigation will benefit from more extensive tabulations of values of the ${\cal C}_{n,p}$ and deeper analytic knowledge of their dependence on both $n$ and $p$. A second requirement will be a much more extensive tabulation of the quantities $\xi_r$ which are needed for the quantities $\Sigma_p^\xi$, using probably the
 Kreminski method \cite{Kreminski}. With these elements in hand, further advances in our understanding of the Riemann hypothesis can be expected.


\begin{thebibliography}{9}
\bibitem{titheath} Titchmarsh, E.C. \& Heath-Brown, D.R., 1986  {\em The Theory of the Riemann Zeta-function}, Oxford, Oxford pp. 254--291.
\bibitem{edw} Edwards, H.M., {\em Riemann's Zeta Function}, Dover, Mineola (2001) pp. 132--134.
\bibitem{plattrud} Platt, D. and Trudgian,T.  2020, The Riemann hypothesis is true up to  $3 \times 10^{12}$,  arXiv:2004.09765v1.
\bibitem{lehmer} Lehmer, D.H. 1988 The sum of like powers of  the zeros of the Riemann zeta function {\em Math. Comp.} {\bf 50} 265--273.
\bibitem{keiper} Keiper, J.B. 1992 Power Series Expansions of Riemann's $\xi$ Function {\em Math. Comp.} {\bf 58} 765-773.
\bibitem{li} Li, X.J. 1997  The positivity of a sequence of numbers and the Riemann hypothesis {\em J. Number Th.} {\bf 65} 325--333.
\bibitem{HAL} McPhedran, R.C., Scott, T.C. and Maignan, A., 2022  Comments on and Extensions to Criteria of Keiper and Li for the Riemann Hypothesis {\em HAL Archive}, hal-03579652v1 and ACM Communications in
Computer Algebra, {\bf 57}, pp. 85-110 (2023).
\bibitem{rmcp2024}  McPhedran, R.C. 2024 Numerical investigations of the Keiper-Li Criterion for the Riemann hypothesis arXiv:2311.06294, pp. 27.
\bibitem {prt} Taylor, P.R. 1945 On the Riemann zeta-function, {\em Q.J.O.}, {\bf 16}, 1-21.
\bibitem{bomblag} Bombieri, E. and Lagarias, J.C. 1999 Complements to Li's Criterion for the Riemann Hypothesis  {\em J. Number Th.} {\bf 77} 274--287.
\bibitem{lagandsuz} Lagarias, J.C. and Suzuki, M. 2006 The Riemann hypothesis for certain integrals of Eisenstein series {\em J. Number Theory} {\bf 118} 98--122.
\bibitem{ki} Ki, H.  2006 Zeros of the constant term in the Chowla-Selberg formula {\em Acta Arithmetica} {\bf 124} 197-204.



\bibitem{jacques3} P\'{o}lya, G. 1926  Bemerkung uber die Integraldarstellung der Riemannschen xi-Funktion
{\em Acta Mathematica},{\bf  48}, 305--317.
\bibitem{pust} Pustyl'nikov, L.D.  1999 On a property of the classical zeta-function associated with the Riemann hypothesis {\em Russian Mathematical Surveys} {\bf 54} 262--263.
\bibitem{pust2} Pustyl'nikov, L.D.  2000 On the asymptotic behaviour of the Taylor series coefficients of $\xi(s)$ {\em Russian Mathematical Surveys} {\bf 55} 349--350.
\bibitem{GORZ} Griffin, M.,  Ono,K.,  Rolen,L. and  Zagier, D. Jensen polynomials for the Riemann zeta function and
other sequences 2019 {\em Proceedings of the National Academy of Sciences}, {\bf 116}, 11103-11110.
\bibitem{kreminski}  Kreminski, R.  2005 \\
http://faculty.colostate-pueblo.edu/rick.kreminski/stieltjesrelated/sigmavalues2005.txt
\bibitem{coffey2005} Coffey, M.W. 2005 Toward verification of the Riemann hypothesis: application of the Li criterion {\em Mathematical Physics, Analysis and Geometry} {\bf 8} 211-255.
\bibitem{titfn}  Titchmarsh, E.C. 1939 {\em  The Theory of Functions}, Oxford, Oxford, Chapter 7.
\bibitem{Kreminski} Kreminski, R. 2003 Newton-Cotes integration for approximating Stieltjes (generalized Euler) coefficients {\em Math. Comp.}, {\bf 72}, 1379-97.
\end{thebibliography}
\end{document}